\documentclass[graybox]{svmult}

\usepackage{type1cm}
\usepackage{makeidx}
\usepackage{graphicx}
\usepackage{multicol}
\usepackage[bottom]{footmisc}

\usepackage{newtxtext}
\usepackage[varvw]{newtxmath}

\usepackage{amsmath}
\usepackage{caption}
\usepackage{subcaption}

\usepackage{cleveref}

\usepackage{todonotes}

\DeclareMathOperator{\diag}{diag}

\makeindex

\begin{document}

\title*{A computational study of algebraic coarse spaces for two-level overlapping additive Schwarz preconditioners}
\titlerunning{Algebraic coarse spaces for two-level Schwarz preconditioners}

\author{Filipe A. C. S. Alves \and Alexander Heinlein \and Hadi Hajibeygi}
\authorrunning{F. A. C. S. Alves \and A. Heinlein \and H. Hajibeygi}

\institute{Filipe A. C. S. Alves \at Delft Institute of Applied Mathematics, Faculty of Electrical Engineering, Mathematics and Computer Science, Delft University of Technology, Delft, The Netherlands. \email{f.a.cumarusilvaalves@tudelft.nl}
\and Alexander Heinlein \at Delft Institute of Applied Mathematics, Faculty of Electrical Engineering, Mathematics and Computer Science, Delft University of Technology, Delft, The Netherlands. \email{a.heinlein@tudelft.nl}
\and Hadi Hajibeygi \at Department of Geoscience and Engineering, Faculty of Civil Engineering and Geosciences, Delft University of Technology, Delft, The Netherlands. \email{h.hajibeygi@tudelft.nl}}

\maketitle

\abstract{The two-level overlapping additive Schwarz method offers a robust and scalable preconditioner for various linear systems resulting from elliptic problems. One of the key to these properties is the construction of the coarse space used to solve a global coupling problem, which traditionally requires information about the underlying discretization. An algebraic formulation of the coarse space reduces the complexity of its assembly. Furthermore, well-chosen coarse basis functions within this space can better represent changes in the problem's properties. Here we introduce an algebraic formulation of the multiscale finite element method (MsFEM) based on the algebraic multiscale solver (AMS) in the context of the two-level Schwarz method. We show how AMS is related to other energy-minimizing coarse spaces. Furthermore, we compare the AMS with other algebraic energy-minimizing spaces: the generalized Dryja--Smith--Widlund (GDSW), and the reduced dimension GDSW (RGDSW).}

\section{Introduction}

Let
\begin{equation} \label{eq:les}
	Au = b
\end{equation}
be the linear equation system arising from the discretization of a boundary value problem on the computational domain $\Omega$. Furthermore, let $\Omega$ be decomposed into $N$ non-overlapping subdomains $\Omega_1,\ldots,\Omega_N$, which we can extend by layers of elements to obtain overlapping subdomains $\Omega_1',\ldots,\Omega_N'$. The corresponding two-level overlapping additive Schwarz preconditioner reads
\begin{equation}
    \label{eq:two_level_precond}
    M_{OAS, 2}^{-1} = \Phi A_0^{-1} \Phi^T + \sum_{i = 1}^N R_i^T A_i^{-1} R_i,
\end{equation}
where $R_i$ are the restriction matrices to the overlapping subdomains $\Omega_i'$, the $A_i = R_i A R_i^T$ are the corresponding local overlapping subdomain matrices, and $A_0 = \Phi^T A \Phi$ is the coarse problem matrix. The key ingredient of the coarse level in the preconditioner above is the choice of prolongation operator $\Phi$, which is composed of the coarse basis functions. In general, the coarse problem is required for numerical scalability of domain decomposition methods; see~\cite{TOSELLI2005DDBOOK} for more details. Additionally, the local overlapping subdomain matrices can be retrieved algebraically; see e.g. \cite{heinlein_fully_2021}.

A first option to define $\Phi$ is by computing a set of finite element shape functions on a coarse triangulation~\cite{TOSELLI2005DDBOOK}. In order to allow for the algebraic construction of $\Phi$, we avoid the use of an additional grid and construct the coarse basis functions based on the nonoverlapping domain decomposition into $\Omega_1,\ldots,\Omega_N$. As a first option, we employ the generalized Dryja--Smith--Widlund (GDSW) coarse space~\cite{DOHRMANN2008GDSW2}. It is composed of energy-minimizing functions of trace functions defined on the interfaces of the nonoverlapping subdomains. The trace values are the restriction of the null space of the global Neumann problem to a decomposition of the interface into vertices, edges, and, in three dimensions, faces. Furthermore, we consider a modification of the GDSW, the reduced dimension GDSW (RGDSW) coarse space introduced in~\cite{DOHRMANN2017RGDSW}. In RGDSW, only vertex functions are considered, resulting in a reduced number of interface components; the interface values are modified by the multiplication of the null space by a partition of unit function.

Related coarse spaces based on the multiscale finite element method (MsFEM)~\cite{HOU1997MSFEM} have been studied in the context of overlapping Schwarz preconditioners~\cite{buck_multiscale_2013}; MsFEM functions are also constructed as energy-minimizing extensions into the interior. They provide good results when applied to media with high coefficient jumps but are generally not algebraic; cf.~\cite{HEINLEIN2018MS} for a more robust variant. A purely algebraic variant of MsFEM coarse spaces was presented in~\cite{LUNATI2009OPMSFV} and later applied in a preconditioning scheme as the two-stage algebraic multiscale solver (TAMS)~\cite{ZHOUTCHELEPI2012TAMS} and the algebraic multiscale solver (AMS) \cite{WANG2014AMS}.

In the following, we establish the relation between AMS and the GDSW and RGDSW coarse basis functions. Our goal is to define a common framework and numerically compare the approaches for some simple heterogeneous model problems. In particular, we will focus on scalar diffusion problems
\begin{alignat}{2}
    \label{eq:diffusion_model}
    - \nabla \cdot \left( \alpha (x) \nabla u (x) \right) &= f(x) \quad \quad && \text{in} \, \Omega \subset \bbbr^d, \\ \nonumber
    u &= u_D (x) \quad \quad && \text{on} \, \partial \Omega,
\end{alignat}
for $d = 2, 3$. Here, the scalar coefficient function $\alpha : \bbbr^d \to \bbbr$ is highly heterogeneous and possibly containing large coefficient jumps.

\section{Energy-minimizing coarse basis functions} \label{sec:coarse_spaces}

Let
\begin{equation*}
    \Gamma = \bigcup_{i=1}^N \partial\Omega_i \setminus \partial\Omega
\end{equation*}
be the interface of the nonoverlapping domain decomposition. We may order the degrees of freedom, such that the system matrix $A$ and the prolongation operator $\Phi$ can be written as
\begin{equation} \label{eq:A_block}
    A = \begin{pmatrix}
        A_{II} & A_{I\Gamma} \\
        A_{\Gamma I} & A_{\Gamma \Gamma}
    \end{pmatrix}, \,
    \Phi = \begin{pmatrix}
        \Phi_I \\
        \Phi_{\Gamma}
    \end{pmatrix},
\end{equation}
where $\Gamma$ indicates the degrees of freedom corresponding to interfaces finite element nodes and $I$ the remaining degrees of freedom; as usual, if the boundary degrees of freedom have not been eliminated, we consider them as interior. The energy-minimizing extensions of interface values $\Phi_{\Gamma}$ are then defined as:
\begin{equation}
    \label{eq:discrete_harmonic_ext}
    \Phi = H_{\Gamma} \Phi_{\Gamma} =
    \begin{pmatrix}
        -A_{II}^{-1} A_{I \Gamma} \\
        I_{\Gamma \Gamma}
    \end{pmatrix} \Phi_{\Gamma},
\end{equation}
where $I_{\Gamma \Gamma}$ is the identity matrix of dimension $| \Gamma | \times | \Gamma |$. \Cref{eq:discrete_harmonic_ext} allows for the description of multiple coarse spaces based on different definitions of the interface values $\Phi_{\Gamma}$.

\subsection{GDSW and RGDSW coarse basis functions} \label{sec:gdsw}

\begin{figure}[t]
	\centering
	\includegraphics[height=0.4\textwidth]{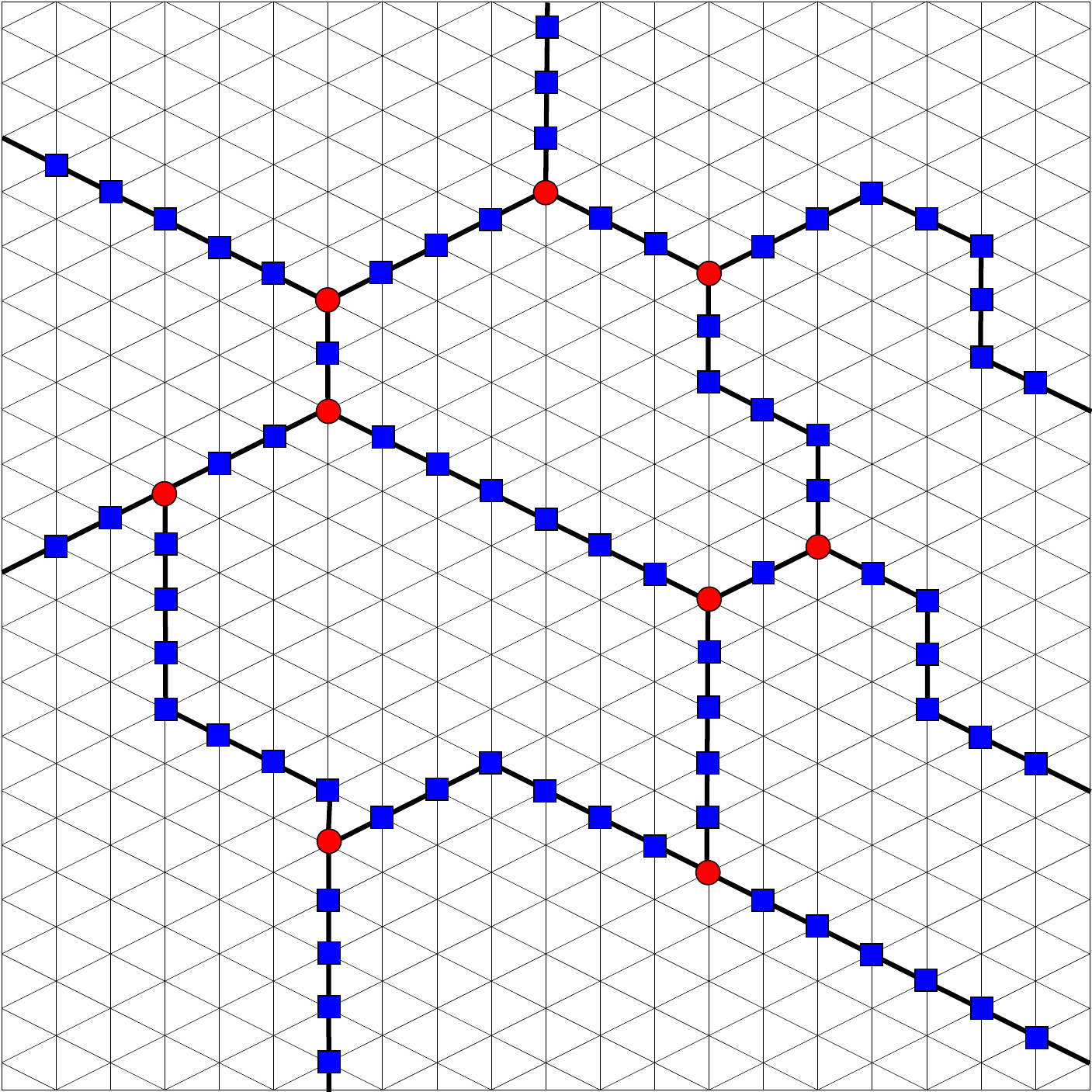}
	\hfill
	\includegraphics[height=0.4\textwidth]{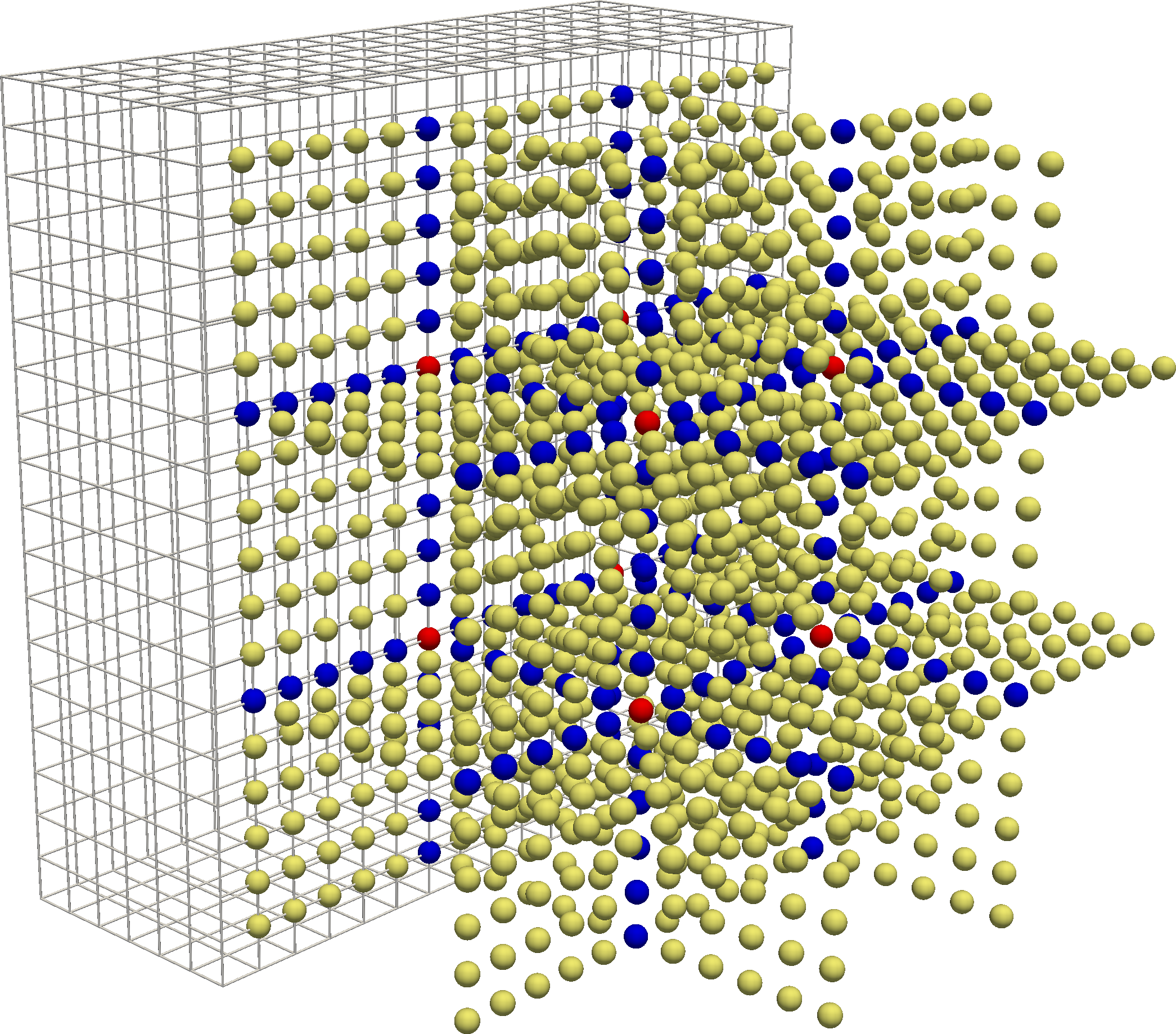}
	\caption{Example of an interface partition of grid elements in two (left) and three (right) dimensions. The vertex, edge and face sub-components are marked in red, blue and yellow, respectively.}
	\label{fig:interface_partition}
\end{figure}

For the GDSW coarse space~\cite{DOHRMANN2008GDSW2}, the interface is decomposed into vertices $V$, edges $E$, and, in three dimensions, faces $F$. An example of decomposition is illustrated in~\cref{fig:interface_partition}. For each interface component, interface values of the corresponding basis function are given by the restriction of the null space of the problem to it. In particular, by the null space of the problem, we mean the null space of the system matrix $A$ with Neumann boundary conditions on the whole boundary $\partial\Omega$. For the RGDSW coarse space~\cite{DOHRMANN2017RGDSW}, the interface is partitioned based on nodal equivalence classes, thus reducing the number of interface components and, consequently, the dimension of the coarse space; in practice, the interface components overlap and mostly correspond to one of the vertices. Moreover, the interface values $\Phi_{\Gamma}$ are modified by multiplication with partition of unit functions. The option 1 described in \cite{DOHRMANN2017RGDSW} can be computed completely algebraically; the interface value in an interface finite element is defined by the inverse of the number of interface components the node belongs to. We will consider this algebraic variant in this paper.

\subsection{AMS coarse basis functions} \label{sec:ams}

AMS~\cite{WANG2014AMS} is a multiscale framework that has been developed for simulations of the Darcy flow in porous media. The AMS coarse basis functions corresponds to an algebraic variant of the MsFEM basis functions~\cite{LUNATI2009OPMSFV}. When applied to a system assembled using a structured grid and a 5-point or 7-point stencil in two or three dimensions, respectively, the resulting basis is the same as for the MsFEM. Otherwise, the AMS gives an approximation of the latter. Some modifications are required in these cases to restrict the support of the coarse basis functions; see \cite{SOUZA2022AMSU}. These changes do not affect the analysis in the context of the two-level Schwarz method.

For the computation of the coarse basis, the grid nodes are again grouped into vertices ($V$), edges ($E$), faces ($F$) and interior ($I$) nodes. This is equivalent to the partition presented in~\cref{fig:interface_partition}. Therefore, let us structure the system~\cref{eq:les} based on those groups of degrees of freedom:
$$
	\begin{pmatrix}
	    A_{II} & A_{IF} & A_{IE} & A_{IV} \\
	    A_{FI} & A_{FF} & A_{FE} & A_{FV} \\
	    A_{EI} & A_{EF} & A_{EE} & A_{EV} \\
	    A_{VI} & A_{VF} & A_{VE} & A_{VV} \\
	\end{pmatrix}
	\begin{pmatrix}
	    u_I \\
	    u_F \\
	    u_E \\
	    u_V
	\end{pmatrix} =
	\begin{pmatrix}
	    b_I \\
	    b_F \\
	    b_E \\
	    b_V
	\end{pmatrix},
$$
where, due to the symmetry of $A$, we have $A_{ij} = A_{ji}^\top$, $i,j \in \left\lbrace I,F,E,V \right\rbrace$.

We require $u$ to be set as the solution of a reduced boundary condition problem on the boundary of the nonoverlapping subdomains. This is referred as the localization assumption \cite{LUNATI2009OPMSFV}. To enforce this assumption in the discrete system, we start by removing the blocks related to connections between internal and interface nodes when solving for edge and face nodes. This results in the system:
\begin{equation}
    \label{eq:ams_reduced_system}
    \begin{pmatrix}
        A_{II} & A_{IF} & A_{IE} & A_{IV} \\
        0 & \tilde{A}_{FF} & A_{FE} & A_{FV} \\
        0 & 0 & \tilde{A}_{EE} & A_{EV} \\
        0 & 0 & 0 & A_0 \\
    \end{pmatrix}
    \begin{pmatrix}
        u_I' \\
        u_F' \\
        u_E' \\
        u_V'
    \end{pmatrix} =
    \begin{pmatrix}
        b_I \\
        b_F \\
        b_E \\
        R_0
    \end{pmatrix}.
\end{equation}
It is important to point out that the solution to the linear system above gives only an approximation $u'$ to the exact solution of \cref{eq:les}. Moreover, $A_0 u_V' = R_0$ corresponds to the coarse-scale problem, which is equivalent to the coarse problem solved in \cref{eq:two_level_precond}. Specifically
\begin{equation*}
    A_0 = \Phi A \Phi^T, \, R_0 = \Phi^T b.
\end{equation*}
where the operator $\Phi$ will be defined later.

The blocks $\tilde{A}_{FF}$ and $\tilde{A}_{EE}$ are modified accordingly to take into account the elimination of the lower triangular part of the system:
\begin{align}
    \label{eq:ams_mass_redist}
    \tilde{A}_{FF} &= A_{FF} + \diag \left( A_{FI} \vec{1}_I \right), \\ \nonumber
    \tilde{A}_{EE} &= A_{EE} + \diag \left( A_{EI} \vec{1}_I \right) + \diag \left( A_{EF} \vec{1}_F \right),
\end{align}
in which $\vec{1}_I$ and $\vec{1}_F$ are vectors of dimension $|I|$ and $|F|$, respectively, that contain only ones. The operation $A_{ij} \vec{1}_j$ corresponds to the sum of the rows of the block $A_{ij}$. In the step above, the diagonal blocks are modified to incorporate the removed off-diagonal blocks in such a way that the energy-minimality of constant functions is maintained. This ensures that corresponding reduced problems generate a partition of unity on the interface, which is an important result for the convergence analysis of the two-level Schwarz preconditioner \cite{TOSELLI2005DDBOOK}.

The system in~\cref{eq:ams_reduced_system} can be solved by backward substitution yielding:
$$
    \begin{pmatrix}
        u_I' & u_F' & u_E' & u_V'
    \end{pmatrix}^T = \Phi u_V' + \mathbf{C},
$$
where $\Phi$ is the AMS prolongation operator written as:
\begin{equation}
    \label{eq:ams_phi}
    \Phi = 
    \begin{pmatrix}
        \Phi_I \\
        \Phi_F \\
        \Phi_E \\
        \Phi_V
    \end{pmatrix} = 
    \begin{pmatrix}
        -A_{II}^{-1} \left( A_{IV} \Phi_V + A_{IE} \Phi_E + A_{IF} \Phi_F \right) \\
        -\tilde{A}_{FF}^{-1} \left( A_{FV} \Phi_V + A_{FE} \Phi_E \right) \\
        -\tilde{A}_{EE}^{-1} A_{EV} \Phi_V \\
        I_{VV}
    \end{pmatrix},
\end{equation}
where $I_{VV}$ is an identity matrix of dimension $|V| \times |V|$.

From~\cref{eq:ams_phi}, it can be seen that the values of the AMS coarse basis functions are computed as energy-minimizing extensions in~\cref{eq:discrete_harmonic_ext}. In particular, $\Phi_I$ in~\cref{eq:ams_phi} can be rewritten as:
$$
    \Phi_I = 
    -A_{II}^{-1} A_{I \Gamma} \Phi_{\Gamma} =
    -A_{II}^{-1}
    \begin{pmatrix}
        A_{IF} & 0 & 0 \\
        0 & A_{IE} & 0 \\
        0 & 0 & A_{IV}
    \end{pmatrix}
    \begin{pmatrix}
        \Phi_F \\
        \Phi_E \\
        \Phi_V
    \end{pmatrix},
$$
which is equivalent to~\cref{eq:discrete_harmonic_ext} when using that $\Gamma = F \cup\, E \cup\, V$.

The AMS coarse space thus only differs from the GDSW and RGDSW coarse spaces in the definition of the interface values of the basis functions. Similar to RGDSW, the AMS basis functions are vertex-centered. However, the values on the interface are computed recursively by solving reduced elliptic problems on the edges and faces; similar recursive computations are also performed in the non-algebraic MsFEM coarse basis functions~\cite{buck_multiscale_2013,HEINLEIN2018MS}. Compared to GDSW coarse spaces and option 1 of RGDSW coarse spaces, where the interface values are constant on each interface component for diffusion problems, the AMS functions are allowed to change according to the contrast in the coefficient function $\alpha$ in~\cref{eq:diffusion_model}.

\section{Numerical examples}

\begin{figure}[t]
	\centering
	\includegraphics[width=0.4\textwidth]{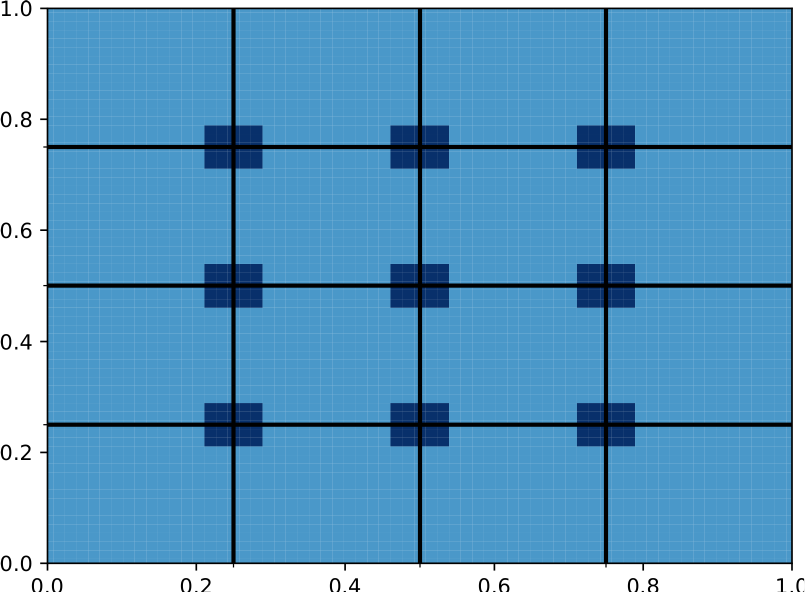}
	\hfill
	\includegraphics[width=0.4\textwidth]{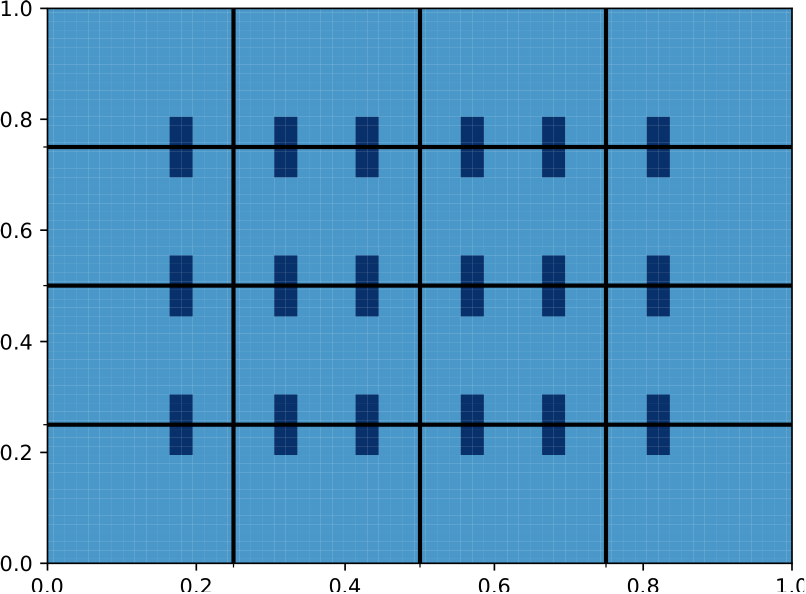}
	\caption{Discontinuous coefficient functions used in the examples. The dark blue regions correspond to the high coefficient inclusions with a value of $\alpha (x) = 10^8$, and $\alpha (x) = 1$ elsewhere. Illustration for $1/H = 4$.}
	\label{fig:coeff_functions}
\end{figure}

In this section, we provide some numerical examples to compare the algebraic coarse spaces discussed in~\cref{sec:coarse_spaces}. For the sake of simplicity, we focus on two-dimensional cases of the model problem~\cref{eq:diffusion_model} on $\Omega = (0,1)^2$, discretized using Q1 finite elements. Furthermore, we only consider zero Dirichlet boundary conditions and a constant source term $f(x) = 1$. The resulting system~\cref{eq:les} is solved using the conjugate gradient (CG) method, preconditioned using the two-level Schwarz preconditioner in~\cref{eq:two_level_precond}, and a relative stopping criterion $|| r^{(k)} ||_2 / || r^{(0)} ||_2 < 10^{-8}$, where $r^{(0)}$ and $r^{(k)}$ are the initial and $k$-th unpreconditioned residuals, respectively. We employ a structured domain decomposition into square subdomains with overlap size $\delta = 2h$ and if not stated otherwise, $H/h = 16$, where $H$ and $h$ are the subdomain and the grid size, respectively. Furthermore, we also present the results using a one-level overlapping additive Schwarz (OAS-1) preconditioner for reference.

The coefficient functions used in the examples are illustrated in~\cref{fig:coeff_functions}. We consider two distributions of high coefficient inclusions: at the coarse nodes (\cref{fig:coeff_functions} (left)), and parallel channels crossing the subdomains interfaces (\cref{fig:coeff_functions} (right)). In both cases, $\alpha (x) = 10^8$ in the dark blue regions and $\alpha (x) = 1$ elsewhere.

\begin{figure}[t]
    \centering
    \includegraphics[width=\textwidth]{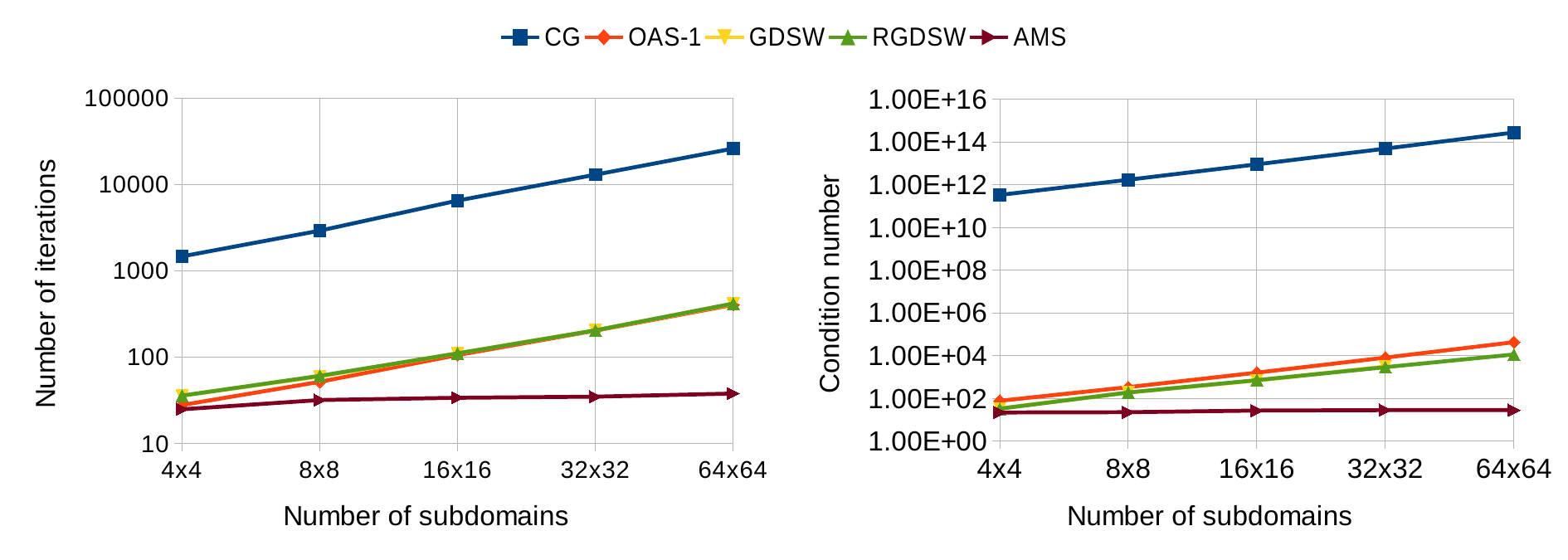}
    \caption{Number of iterations (left) and condition number estimate (right) versus the number of subdomains for the coefficient function in~\cref{fig:coeff_functions} (left).}
    \label{fig:ex2_results}
\end{figure}

The simulation results for the coefficient function in~\cref{fig:coeff_functions} (left) are presented in~\cref{fig:ex2_results}. In this scenario, the GDSW and the RGDSW coarse spaces are not robust with respect to the coefficient: the condition is in the order of the contrast and the iteration counts increase with the number of subdomains. On the other hand, the AMS coarse space is robust and scalable. This is in line with previous observations for coarse spaces based on MsFEM \cite{HEINLEIN2018MS}. This family of coarse spaces remain robust as long as the high coefficient inclusions are connected to the coarse nodes.

\begin{figure}[t]
    \centering
    \includegraphics[width=\textwidth]{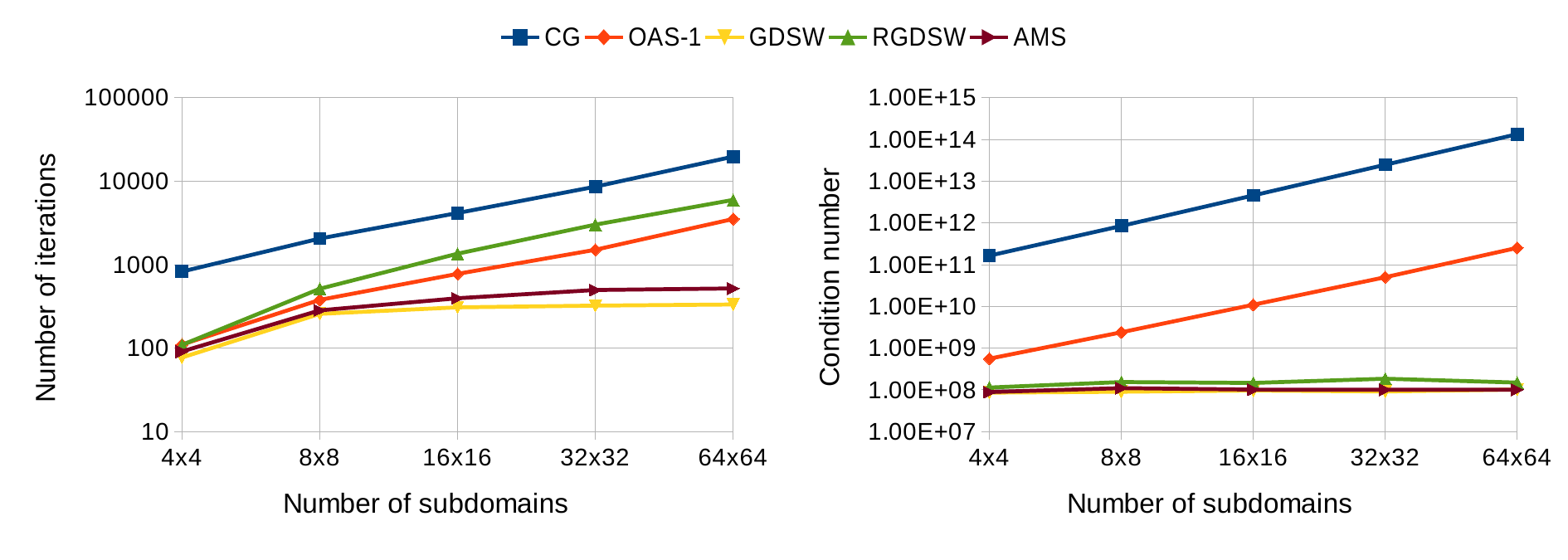}
    \caption{Number of iterations (left) and condition number estimate (right) versus the number of subdomains for the coefficient function in~\cref{fig:coeff_functions} (right).}
    \label{fig:ex3_results}
\end{figure}

In~\cref{fig:ex3_results}, the results for the coefficient function in~\cref{fig:coeff_functions} (right) are shown. Here, all three coarse spaces yield a condition number in the order of the contrast. As is well-known, the coefficient contrast across the edge cannot be solved by the vertex-centered functions from the RGDSW and AMS coarse spaces, while the GDSW basis functions cannot deal with multiple jumps across the edges. However, we notice that GDSW and AMS perform better than RGDSW in terms of the number of iterations, and even though iteration counts are high, they seem to stabilize with increasing numbers of subdomains.

\begin{figure}[!ht]
    \centering
    \includegraphics[width=\textwidth]{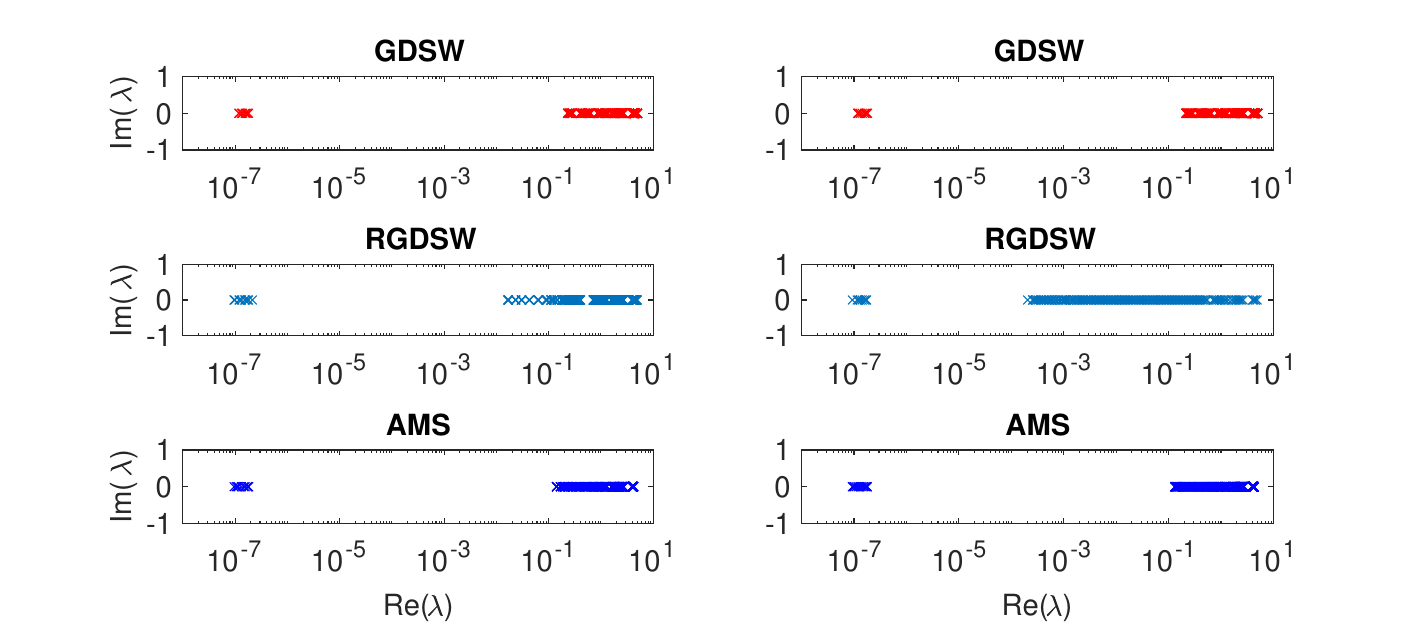}
    \caption{Eigenvalues for each coarse space for the coefficient function in~\cref{fig:coeff_functions} (right) and different subdomain sizes: \textbf{left}: $1/H = 8$; \textbf{right}: $1/H = 64$.}
    \label{fig:ex3_eig_vals}
\end{figure}

In order to analyze the difference in the convergence behavior of RGDSW versus GDSW and AMS coarse spaces, we investigate the spectrum of the preconditioned system matrix for the three coarse spaces in~\cref{fig:ex3_eig_vals}. We observe that the low eigenvalues are clustered around $10^{-7}$, which remains the same when increasing the number of subdomains. Thus the condition number also remains roughly constant. We observe that, while the clustering of the spectrum remains the same for the AMS and GDSW coarse spaces, the spectrum for the RGDSW coarse space is stretched significantly when increasing the number of subdomains. This could explain the different behavior in the convergence for increasing numbers of subdomains.

\begin{acknowledgement}
This study was funded by Energi Simulation and the Khalifa University through Project MIRE (Multiscale Integrated Reservoir Earth).
\end{acknowledgement}

\bibliographystyle{spmpsci}
\bibliography{bibliography}

\begin{thebibliography}{10}
\providecommand{\url}[1]{{#1}}
\providecommand{\urlprefix}{URL }
\expandafter\ifx\csname urlstyle\endcsname\relax
  \providecommand{\doi}[1]{DOI~\discretionary{}{}{}#1}\else
  \providecommand{\doi}{DOI~\discretionary{}{}{}\begingroup \urlstyle{rm}\Url}\fi

\bibitem{buck_multiscale_2013}
Buck, M., Iliev, O., Andrä, H.: Multiscale finite element coarse spaces for the application to linear elasticity.
\newblock Open Mathematics \textbf{11}(4), 680--701 (2013).
\newblock Publisher: De Gruyter Open Access

\bibitem{DOHRMANN2008GDSW2}
Dohrmann, C.R., Klawonn, A., Widlund, O.B.: {Domain Decomposition for Less Regular Subdomains: Overlapping Schwarz in Two Dimensions}.
\newblock SIAM J. Numer. Anal. \textbf{46}, 2153--2168 (2008)

\bibitem{DOHRMANN2017RGDSW}
Dohrmann, C.R., Widlund, O.B.: {On the Design of Small Coarse Spaces for Domain Decomposition Algorithms}.
\newblock SIAM J. Sci. Comput. \textbf{39}, A1466--A1488 (2017)

\bibitem{heinlein_fully_2021}
Heinlein, A., Hochmuth, C., Klawonn, A.: Fully algebraic two-level overlapping {Schwarz} preconditioners for elasticity problems.
\newblock In: Numerical mathematics and advanced applications—{ENUMATH} 2019, \emph{Lect. {Notes} {Comput}. {Sci}. {Eng}.}, vol. 139, pp. 531--539. Springer, Cham (2021)

\bibitem{HEINLEIN2018MS}
Heinlein, A., Klawonn, A., Knepper, J., Rheinbach, O.: {Multiscale coarse spaces for overlapping Schwarz methods based on the ACMS space in 2D}.
\newblock {Electron. Trans. Numer. Anal.} \textbf{48}, 156--182 (2018)

\bibitem{HOU1997MSFEM}
Hou, T.Y., Wu, X.H.: {A Multiscale Finite Element Method for Elliptic Problems in Composite Materials and Porous Media}.
\newblock J. Comput. Phys. \textbf{134}, 169--189 (1997)

\bibitem{LUNATI2009OPMSFV}
Lunati, I., Lee, S.H.: {An Operator Formulation of the Multiscale Finite-Volume Method with Correction Function}.
\newblock Multiscale Model. Simul. \textbf{8}, 96--109 (2009)

\bibitem{SOUZA2022AMSU}
Souza, A.C.R., Carvalho, D.K.E., Santos, J.C.A., Willmersdorf, R.B., Lyra, P.R.M., Edwards, M.G.: {An algebraic multiscale solver for the simulation of two-phase flow in heterogeneous and anisotropic porous media using general unstructured grids (AMS-U)}.
\newblock Appl. Math. Modell. \textbf{103}, 792--823 (2022)

\bibitem{TOSELLI2005DDBOOK}
Toselli, A., Widlund, O.B.: {Domain Decomposition Methods — Algorithms and Theory}, vol.~34.
\newblock Springer Berlin Heidelberg (2005)

\bibitem{WANG2014AMS}
Wang, Y., Hajibeygi, H., Tchelepi, H.A.: {Algebraic multiscale solver for flow in heterogeneous porous media}.
\newblock J. Comput. Phys. \textbf{259}, 284--303 (2014)

\bibitem{ZHOUTCHELEPI2012TAMS}
Zhou, H., Tchelepi, H.A.: {Two-Stage Algebraic Multiscale Linear Solver for Highly Heterogeneous Reservoir Models}.
\newblock SPE Journal \textbf{17}(02), 523--539 (2012)

\end{thebibliography}

\end{document}